\input amstex\documentstyle{amsppt}  
\pagewidth{12.5cm}\pageheight{19cm}\magnification\magstep1
\topmatter
\title Parabolic character sheaves, III\endtitle
\author G. Lusztig\endauthor
\address{Department of Mathematics, M.I.T., Cambridge, MA 02139}\endaddress
\thanks{Supported in part by the National Science Foundation}\endthanks
\endtopmatter   
\document

\define\hK{\hat K}

\define\po{\text{\rm pos}}

\define\si{\sim}

\define\sqc{\sqcup}

\define\qua{\quad}

\define\part{\partial}
\define\em{\emptyset}

\define\iy{\infty}
\define\m{\mapsto}
\define\do{\dots}

\define\bsl{\backslash}

\define\sub{\subset}    

\define\T{\times}
\define\ti{\tilde}
\define\nl{\newline}
\redefine\i{^{-1}}

\define\ov{\overline}

\define\bbq{\bar{\QQ}_l}

\define\Ad{\text{\rm Ad}}

\define\ind{\text{\rm ind}}

\redefine\b{\beta}
\redefine\c{\chi}
\define\g{\gamma}
\redefine\d{\delta}
\define\e{\epsilon}

\redefine\o{\omega}

\redefine\t{\tau}

\define\vt{\vartheta}

\redefine\D{\Delta}

\define\kk{\bold k}
\redefine\ll{\bold l}

\redefine\tt{\bold t}

\define\GG{\bold G}

\define\NN{\bold N}

\define\QQ{\bold Q}

\define\WW{\bold W}
\define\ZZ{\bold Z}

\define\ca{\Cal A}
\define\cb{\Cal B}
\define\cc{\Cal C}

\define\cl{\Cal L}

\define\co{\Cal O}
\define\cp{\Cal P}

\define\ct{\Cal T}

\define\tK{\ti K}

\define\sh{\sharp}

\define\sps{\supset}

\define\sneq{\subsetneqq}
\define\BE{B}
\define\LI{L1}
\define\LII{L2}
\define\LIII{L3}
\define\LIV{L4}

\head 1. A decomposition of $\GG^\d/U_P$\endhead
\subhead 1.1\endsubhead
Let $\kk$ be an algebraically closed field. In this paper an algebraic variety (or algebraic group) is always 
assumed to be over $\kk$. More generally we will consider ind-varieties. An ind-variety is an increasing union 
$X=\cup_{i\in\NN}X_i$ of algebraic varieties $X_i$ under closed imbeddings $j_i:X_i@>>>X_{i+1}$. (The sequence 
$(X_i)$ can be replaced by a sequence $(X'_i)$ such that each $X_i$ is a closed subvariety of some $X'_j$ and each
$X'_i$ is a closed subvariety of some $X_j$. This change leads to the same ind-variety.) An algebraic subvariety 
of $X$ is by definition an algebraic subvariety of $X_i$ for some $i$. We fix a prime number $l$ invertible in 
$\kk$. For $X$ as above we can define the category of perverse sheaves on $X$. Its objects are collections 
$K=(K_i)_{i\ge i_0}$ where $K_i$ is a $\bbq$-perverse sheaf on $X_i$ for $i\ge i_0$, such that 
$K_{i+1}=j_{i!}K_i$ for any $i\ge i_0$. (The support of $K$ is by the definition the support of $K_i$ for 
$i\ge i_0$. $K$ is said to be simple if $K_i$ is simple for $i\ge i_0$.)

The main purpose of this paper is to define a class of simple perverse sheaves (called character sheaves) on 
certain ind-varieties associated to a loop group. This has applications to a geometric construction of certain 
affine Hecke algebras with unequal parameters (an affine analogue of a construction in \cite{\LIV}), as will be 
shown elsewhere.

Let $\ll=\kk((\e))$ where $\e$ is an indeterminate. Let $G$ be a connected semisimple almost simple algebraic 
group. Let $\GG=G(\ll)$ and let $\GG'$ be the derived group of $\GG$. A subgroup of $\GG'$ is said to be parahoric
if it contains some Iwahori subgroup and is not equal to $\GG'$.
Let $P$ be a parahoric subgroup of $\GG'$. Note that $P$ can be naturally regarded as a proalgebraic group; it has
a canonical normal subgroup $U_P$ which is prounipotent and the quotient $P/U_P$ is naturally a connected
reductive algebraic group. Now $\GG'/U_P$ can be viewed naturally as an ind-variety; indeed $\GG'/U_P$ is fibred 
over $\GG'/P$ (which is known to be an ind-variety) with fibres isomorphic with the algebraic variety $P/U_P$.  
Note that $P$ acts on $\GG'/U_P$ by $p:gU_P\m pgp\i U_P$. In this paper we define a class of simple perverse
sheaves (character sheaves) on $\GG'/U_P$; these are are certain $P$-equivariant simple perverse sheaves on 
$\GG'/U_P$ supported by $P$-stable algebraic subvarieties (on which $P$ acts algebraically through a quotient 
which is an algebraic group). The key to the definition of character sheaves on $\GG'/U_P$ is a decomposition of 
$\GG'/U_P$ into $P$-stable smooth algebraic subvarieties indexed by a subset of the affine Weyl group. This 
decomposition is the affine analogue of a decomposition which appeared in \cite{\LI}; the character sheaves that
we define are affine analogues of the parabolic character sheaves of \cite{\LI} (see also \cite{\LII}).

We will actually consider character sheaves on an ind-variety slightly more general that $\GG'/U_P$, namely 
$\GG^\d/U_P$ where $\d$ is an element of the finite abelian group $D=\GG/\GG'$ and $\GG^\d$ is the inverse image 
of $\d$ under the obvious map $\GG@>>>D$.

Now many of the results in \cite{\LI, \S3} remain valid (with the same proof) if $G,G^1$ are replaced by $\GG'$,
$\GG^\d$, "Borel" is replaced by "Iwahori" and "parabolic" is replaced by "parahoric". We will sometime refer to a
result in \cite{\LI} such as \cite{\LI, 3.2} by \cite{\LI,"3.2"} with the understanding that the replacements in 
the previous sentence are made.

{\it Notation.} If $P,Q$ are parahoric subgroups of $\GG'$ then so is $P^Q=(P\cap Q)U_P$ and we have 
$U_{(P^Q)}=U_P(P\cap U_Q)$. We say that $P,Q$ are in good position if $P^Q=P$ or equivalently $Q^P=Q$.

\subhead 1.2\endsubhead
Let $\cb$ be the set of all Iwahori subgroups of $\GG$ (they are actually contained in $\GG'$). Let $\WW'$ be the
set of $\GG'$-orbits on $\cb\T\cb$ ($\GG'$ acts by conjugation on both factors). For $B,B'\in\cb$ we write 
$\po(B,B')=w$ if the $\GG'$-orbit of $(B,B')$ is $w$. If $w\in\WW$ and $\po(B,B')=w$ then $B/(B\cap B')$ is an 
algebraic variety of dimension say $l(w)$. Let $I=\{w\in\WW';l(w)=1\}$. There is a well defined group 
structure on $\WW'$ in which the elements of $I$ have order $2$ and in which the following property holds: if 
$\po(B,B')=w_1,\po(B',B'')=w_2$ and $l(w_1w_2)=l(w_1)+l(w_2)$ then $\po(B,B'')=w_1w_2$. This makes $\WW'$ into an
infinite Coxeter group (an affine Weyl group) in which the length function is $w\m l(w)$. For $J\sneq I$ let 
$\WW'_J$ be the (finite) subgroup of $\WW'$ generated by $J$. For $J\sneq I$, $K\sneq I$ let ${}^K\WW',\WW'{}^J$,
${}^K\WW'{}^J$ be as in \cite{\LI, 2.1} (with $W$ replaced by $\WW'$).

If $P$ is a parahoric subgroup of $\GG'$, the set of all $w\in\WW'$ such that $w=\po(B,B')$ for some $B,B'\in\cb$,
$B\sub P,B'\sub P$ is of the form $\WW'_J$ for a well defined subset $J\sneq I$. We then say that $P$ is of type 
$J$. For $J\sneq I$ let $\cp_J$ be the set of all parahoric subgroups of type $J$. Then $\GG'$ acts transitively
(by conjugation) on $\cp_J$. For $P\in\cp_J$, $Q\in\cp_K$, there is a well defined element 
$u=\po(P,Q)\in{}^J\WW'{}^K$ such that 
$\po(B,B')\ge u$ (standard partial order on $\WW'$) for any $B,B'\in\cb$, $B\sub P$, $B'\sub Q$ and 
$\po(B_1,B'_1)=u$ for some $B_1,B'_1\in\cb$, $B_1\sub P$, $B'_1\sub Q$; we then have $B_1\sub P^Q$, 
$B'_1\sub Q^P$. We have $P^Q\in\cp_{J\cap\Ad(u)K}$. Also, $(P,Q)\m u$ defines a bijection between the set of 
$\GG'$-orbits on $\cp_J\T\cp_K$ and ${}^J\WW'{}^K$.

\subhead 1.3\endsubhead
Throughout this paper we fix $\d\in D$. There is a unique automorphism of $\WW'$ (denoted again by $\d$) such that
$\d(I)=I$ and such that for any $J\sub I$, any $P\in\cp_J$ and any $g\in\GG^\d$ we have $gPg\i\in\cp_{\d(J)}$. In 
this subsection we recall some results of B\'edard \cite{\BE} (see also \cite{\LI, \S2}). For any $J\sneq I$ let 
$\ct(J,\d)$ be the set of all sequences $\tt=(J_n,w_n)_{n\ge0}$ where $J=J_0\sps J_1\sps J_2\sps\do$ and 
$w_0,w_1,\do$ are elements of $\WW'$ such that \cite{\LI, 2.2(a)-(d)} hold (with $\WW',\d$ instead of $W,\e$). For
$(J_n,w_n)_{n\ge0}\in\ct(J,\d)$ we have $J_n=J_{n+1}=\do=J_\iy$, $w_n=w_{n+1}=\do=w_\iy$ for $n$ large and 
$\Ad(w_\iy)J_\iy=\d(J_\iy)$. Moreover, $(J_n,w_n)_{n\ge0}\m w_\iy$ is a bijection 
$$\ct(J,\d)@>\si>>{}^{\d(J)}\WW'.\tag a$$ 

\subhead 1.4\endsubhead
We fix $J\sneq I$. Let 
$$Z_{J,\d}=\{(P,P',gU_P);P\in\cp_J,P'\in\cp_{\d(J)},gU_P\in\GG^\d/U_P,gPg\i=P'\}.$$
To any $(P,P',gU_P)\in Z_{J,\d}$ we associate a sequence $(J_n,w_n)_{n\ge0}$ with $J_n\sneq I$, $w_n\in\WW'$, and
a sequence $(P^n,P'{}^n,gU_{P^n})\in Z_{J_n,\d}$. We set
$$P^0=P,P'{}^0=P',\qua J_0=J,\qua w_0=\po(P'{}^0,P^0).$$
Assume that $n\ge1$, that $P^m,P'{}^m,J_m,w_m$ are already defined for $m<n$ and that $w_m=\po(P'{}^m,P^m)$ for 
$m<n$, $(P^m,P'{}^m,gU_{P^m})\in Z_{J_m,\d}$ for $m<n$. Let
$$J_n=J_{n-1}\cap\d\i\Ad(w_{n-1})(J_{n-1}),$$
$$P^n=g\i((P'{}^{n-1})^{P^{n-1}})g\in\cp_{J_n},\qua P'{}^n=(P'{}^{n-1})^{P^{n-1}}\in\cp_{\d(J_n)},$$
$$w_n=\po(P'{}^n,P^n)\in{}^{\d(J_n)}W^{J_n}.$$
Note that $gU_{P^n}$ is well defined by $gU_{P^{n-1}}$ since $P^n\sub P^{n-1}$ hence $U_{P^{n-1}}\sub U_{P^n}$. 
This completes the inductive definition. Note that $P'=P'{}^0\sps P'{}^1\sps\do$ and $P=P^0\sps P^1\sps\do$ hence
$(P^n,P'{}^n)$ is independent of $n$ for large $n$. From \cite{\LI, "3.2(c)"} (with $P,P',Z$ replaced by 
$P^{n-1},P'{}^{n-1},P^n$) we see that $w_n\in w_{n-1}\WW'_{J_{n-1}}$ for $n\ge1$. Thus, 
$(J_n,w_n)_{n\ge0}\in\ct(J,\d)$. We write $(J_n,w_n)_{n\ge0}=\b'(P,P',gU_P)$. For $\tt\in\ct(J,\d)$ we set
$${}^\tt Z_{J,\d}=\{(P,P',gU_P)\in Z_{J,\d};\b'(P,P',gU_P)=\tt\}.$$
The $\GG'$-action on $Z_{J,\d}$ given by 
$$h:(P,P',gU_P)\m(hPh\i,hP'h\i,hgh\i U_{hPh\i})$$ 
preserves the subset ${}^\tt Z_{J,\d}$. Clearly, $({}^\tt Z_{J,\d})_{\tt\in\ct(J,\d)}$ is a partition of 
$Z_{J,\d}$. 

\subhead 1.5\endsubhead
Let $J\sneq I$. Let $P\in\cp_J$. The inclusion $a:\GG^\d/U_P@>>>Z_{J,\d}$, $a(gU_P)=(P,gPg\i,gU_P)$ is 
$P$-equivariant where $P$ acts on $\GG^\d/U_P$ by $p:gU_P\m pgp\i U_P$ and on $Z_{J,\d}$ by restriction of the
$\GG'$-action. For any $\tt\in\ct(J,\d)$ we set ${}^\tt\GG^\d/U_P=a\i({}^\tt Z_{J,\d})$. Clearly,
$({}^\tt\GG^\d/U_P)_{\tt\in\ct(J,\d)}$ is a partition of $\GG^\d/U_P$ into $P$-stable subsets.

For example, if $J=\em$, $P\in\cb$ and we identify $\ct(\em,\d)=\WW'$ as in 1.3(a), then for $w\in\WW'$ we have 
${}^w\GG^\d/U_B=\{gU_B\in\GG^\d/U_B;\po(gBg\i,B)=w\}$.

\subhead 1.6\endsubhead
Let $\tt\in Z_{J,\d}$. For any $r\ge0$ let $\tt_r=(J_n,w_n)_{n\ge r}\in\ct(J_r,\d)$. Consider the
sequence of $\GG'$-equivariant maps
$${}^{\tt_0}Z_{J_0,\d}@>\vt_0>>{}^{\tt_1}Z_{J_1,\d}@>\vt_1>>{}^{\tt_2}Z_{J_2,\d}@>\vt_2>>\do\tag a$$
where for $r\ge0$ and any $(P,P',gU_P)\in{}^{\tt_r}Z_{J_r,\d}$ we set
$$\vt_r(Q,Q',gU_Q)\m(g\i(Q'{}^Q)g,Q'{}^Q,gU_{g\i(Q'{}^Q)g})\in{}^{\tt_{r+1}}Z_{J_{r+1},\d}.$$
Note that for sufficiently large $r$, $\vt_r$ is the identity map. (Recall that
$J_r=J_{r+1}=\do=J_\iy$, $w_r=w_{r+1}=\do=w_\iy$ for $r$ large and $\Ad(w_\iy)(J_\iy)=\d(J_\iy)$.) Note that 

(b) {\it  each fibre of $\vt_r$ is isomorphic to the affine space $(U_Q\cap U_{Q'})\bsl(U_Q\cap Q')$ (a closed 
subset of the algebraic variety $U_{Q'}\bsl Q'$) where}

$(Q,Q')\in\cp_{J_r}\T\cp_{\d(J_r)},\po(Q',Q)=w_r$.
\nl
(See \cite{\LI, "3.12(b)"}.) Also,

(c) {\it the map induced by $\vt_r$ from the set of $\GG'$-orbits on ${}^{\tt_r}Z_{J_r,\d}$ to the set of 
$\GG'$-orbits on ${}^{\tt_{r+1}}Z_{J_{r+1},\d}$ is a bijection.}
\nl
(See \cite{\LI, "3.12(c)"}.)

Let $P\in\cp_J$. For any $r\ge0$ we set $\D_r=\{(Q,Q',gU_Q)\in{}^{\tt_r}Z_{J_r,\d};Q\sub P\}$. Now $P$ acts on 
$\D_r$ by restriction of the $\GG'$-action on ${}^{\tt_r}Z_{J_r,\d}$. We have a cartesian diagram
$$\CD
\D_r@>c_r>>{}^{\tt_r}Z_{J_r,\d}\\
@Vb_rVV       @V\vt_rVV\\ 
\D_{r+1}@>c_{r+1}>>{}^{\tt_{r+1}}Z_{J_{r+1},\d}
\endCD$$
where $c_r,c_{r+1}$ are the obvious inclusions and $b_r$ is the restriction of $\vt_r$. Using this and (b) we see
that

(d) {\it each fibre of the $P$-equivariant map $b_r:\D_r@>>>\D_{r+1}$ is isomorphic to the affine space
$(U_Q\cap U_{Q'})\bsl(U_Q\cap Q')$ where $(Q,Q')\in\cp_{J_r}\T\cp_{\d(J_r)},\po(Q',Q)=w_r$.}
\nl
We show:

(e) {\it the map induced by $b_r$ from the set of $P$-orbits on $\D_r$ to the set of $P$-orbits on $\D_{r+1}$ is a
bijection.}
\nl
The cartesian map above induces a commutative diagram
$$\CD
(P-\text{orbits on }\D_r)@>>>(\GG'-\text{orbits on }{}^{\tt_r}Z_{J_r,\d})\\
@VVV       @VVV\\ 
(P-\text{orbits on }\D_{r+1})@>>>(\GG'-\text{orbits on }{}^{\tt_{r+1}}Z_{J_{r+1},\d})
\endCD$$
The horizontal maps are clearly bijections. The right vertical map is a bijection by (c). It follows that the left
vertical map is a bijection. This proves (e).

If $\kk$ is replaced by a finite field $F_q$ with $q$ elements and if the Frobenius map acts trivially on
$\WW',\d$, then ${}^\tt\GG^\d/U_P$ becomes a finite set with cardinal equal to $\sh(P/U_P)q^{l(w)}$ where 
$w\in{}^{\d(J)}\WW'$ corresponds to $\tt$ under 1.3(a). (This is seen by an argument similar to that in 
\cite{\LII, 8.20}.)

\subhead 1.7\endsubhead
We now make a digression. Let $L$ be a group and let $E$ be a set with a free transitive left $L$-action 
$(l,e)\m le$, and a free transitive right $L$-action $(e,l)\m el$ such that $(le)l'=l(el')$ for 
$l,l'\in L,e\in E$. For any $e\in E$ we define a map $\t_e:L@>>>L$ by $\t_e(l)e=el$. Note that $\t_e$ is a group 
automorphism of $L$. Assume further that $E_0:=\{e\in E;\t_e^N=1\text{ for some }N\ge1\}$ is nonempty. Let 
$e_0\in E_0$ and let $d\ge1$ be such that $\t_{e_0}^d=1$. We consider the group $L\g_d$ (semidirect product) where
$\g_d$ is the cyclic group of order $d$ with generator $\o$ and $\o l=\t_{e_0}(l)\o$ for any $l\in L$. Define a 
bijection $f:E@>\si>>L\o$ ($L\o$ is an $L$-coset in $L\g_d$) by $le_0\m l\o$. For any $l\in L,e\in E$ we have 
$f(lel\i)=lf(e)l\i$ (in the right side we use the group structure of $L\g_d$). Thus, the set $E$ with the 
$L$-action $l:e\m lel\i$ is isomorphic to the coset $L\o$ in $L\g_d$ with the $L$-action given by $L$-conjugacy in
the group structure of $L\g_d$.

\subhead 1.8\endsubhead 
We return to the setup in 1.6. For $r\ge0$ sufficiently large, $\tt_r$ is independent of $r$; we denote it by
$\tt_\iy$. We have $\tt_\iy=(J'_n,w'_n)_{n\ge0}$ where $J'_n=J_\iy$, $w'_n=w_\iy$ for all $n\ge0$ and 
$\Ad(w_\iy)(J_\iy)=\d(J_\iy)$. It follows that 
$${}^{\tt_\iy}Z_{J_\iy,\d}=\{(Q,Q',gU_Q);(Q,Q')\in\co,gU_Q\in\GG^\d/U_Q,gQg\i=Q'\}$$
where
$$\co=\{(Q,Q')\in\cp_{J_\iy}\T\cp_{\d(J_\iy)};\po(Q',Q)=w_\iy\}.$$
Note that any $Q,Q'$ with $(Q,Q')\in\co$ are in good position. Hence 
$$\D:=\{(Q,Q',gU_Q);(Q,Q')\in\co,gU_Q\in\GG^\d/U_Q,gQg\i=Q',Q\sub P\}.$$
We have maps $\D@>e>>\D'@>e'>>\D''$ where $\D'=\{(Q,Q')\in\co;Q\sub P\}$, $\D''=\{Q\in\cp_{J_\iy};Q\sub P\}$,
$e(Q,Q',gU_Q)=(Q,Q')$, $e'(Q,Q')=Q$. We can regard $\D,\D',\D''$ as smooth algebraic varieties so that $\D''$ is 
isomorphic to a partial flag manifold $P/Q$ (where $Q\in\D''$) of $P/U_P$, $e'$ is an affine space bundle and $e$
is a fibration with fibres isomorphic to $Q/U_Q$ (where $Q\in\D''$). We can regard $\D_0,\D_1,\D_2,\do$ as smooth
algebraic varieties such that $b_r:\D_r@>>>\D_{r+1}$ is an affine space bundle for any $r\ge0$ (note that $\D_r=\D$ for $r$ large is a smooth algebraic variety). We can identify ${}^\tt\GG^\d/U_P$ with
$\D_0$ by $gU_P\m(P,gPg\i,gU_P)$. It follows that ${}^\tt\GG^\d/U_P=\D_0$ is naturally a smooth algebraic variety,
an iterated affine space bundle over $\D$ via the map 
$$b:=\do b_2b_1b_0:{}^\tt\GG^\d/U_P@>>>\D.$$
In fact, it is an algebraic subvariety of the ind-variety $\GG^\d/U_P$.

Now the $P$-action on $\D$ induces a transitive $P$-action on $\D'$. Hence if we fix $(Q,Q')\in\D'$ and we set
$E=\{gU_Q\in\GG^\d/U_Q,gQg\i=Q'\}$, we have a $P$-equivariant isomorphism 
$$\c:P\T_{Q\cap Q'}E@>\si>>\D,(p,gU_Q)\m(pQp\i,pQ'p\i,pgp\i U_{pQp\i}).$$
Let $L_\tt=(Q\cap Q')/(U_Q\cap U_{Q'})$. Since $Q,Q'$ are in good position, the obvious homomorphisms 
$Q/U_Q@<i<<L_\tt@>i'>>Q'/U_{Q'}$ are isomorphisms. 

Note that the $(Q\cap Q')$-action on $E$ given by $h:gU_Q\m hgh\i U_Q$ factors through an $L_\tt$-action (called 
"conjugation"). On $E$ we have a free transitive right $L_\tt$-action given by $l:gU_Q\m\m gU_Q\cdot l=gqU_Q$ 
(where $qU_Q=i(l)$) and a free transitive left $L_\tt$-action given by $l:gU_Q\m l\cdot gU_Q=q'gU_Q$ (where 
$q'U_{Q'}=i'(l)$), so that the "conjugation" $L_\tt$-action above is $l:gU_Q\m l\cdot gU_Q\cdot l\i$. The left and
right $L_\tt$-actions on $E$ satisfy the hypotheses of 1.7 with $L=L_\tt$. Hence, by 1.7, the conjugation 
$L_\tt$-action on $E$ is isomorphic to the $L_\tt$-action on the connected component $C_\tt=L_\tt\o$ of a 
semidirect product $L_\tt\g_d$ (an algebraic group with identity component $L_\tt$ and with cyclic group of 
components $\g_d$ of order $d$ with generator $\o$).

Thus $\c$ can be viewed as an isomorphism $\c':P\T_{Q\cap Q'}C_\tt@>\si>>\D$ where $Q\cap Q'$ acts on $C_\tt$ via
its quotient $L_\tt$, by conjugation. 

This isomorphism induces a bijection between the set of $P$-orbits on $\D$ and the set of $P$-orbits on 
$P\T_{Q\cap Q'}C_\tt$ which can be naturally identified with the set of $L_\tt$-orbits (for conjugation) on 
$C_\tt$. Composing this with the bijections in 1.6(e) we obtain a bijection between the set of $P$-orbits on 
${}^\tt\GG^\d/U_P=\D_0$ and the set $L_\tt\bsl C_\tt$ of $L_\tt$-orbits (for conjugation) on $C_\tt$. Putting 
together these bijections (for various $\tt$) we obtain a bijection between the set of $P$-orbits on $\GG^\d/U_P$
and the set $\sqc_{\tt\in\ct(J,\d)}L_\tt\bsl C_\tt$.

\subhead 1.9\endsubhead
Let $V$ be a $2$-dimensional vector space over $\ll=\kk((\e))$ with a fixed volume element $\o$. Let 
$\ca=\kk[[\e]]$. Let $X'$ be the set of all $\ca$-lattice in $V$. Let $X$ be the set of all $\cl\in X'$ such that
for some $\ca$-basis $e_1,e_2$ of $\cl$ we have $e_1\wedge e_2=\o$. Then $\GG=SL(V)$ acts transitively on $X$. Let
$\cl\in X$ and let $P=\{g\in\GG;g\cl=\cl\}$. Note that $\GG=\GG'$ is as in 1.1 and $P$ is a parahoric subgroup of
$\GG$. In our case the set $\GG'/U_P$ may be identified with the set $Y$ of all pairs $(\cl',g)$ where $\cl'\in X$
and $g$ is an isomorphism of $\kk$-vector spaces $\cl/\e\cl@>\si>>\cl'/\e\cl'$ preserving the volume elements 
induced by $\o$. We have a partition $Y=\sqc_{n\in\NN}Y_n$ where 
$Y_n=\{(\cl',g)\in Y;\dim_\kk\cl/(\cl\cap\cl')=n\}$. If $n>0$ then $Y_n$ has a further partition 
$Y_n=Y'_n\sqc Y''_n$ defined as folows. Let $(\cl',g)\in Y_n$. Let $\cl_1$ be the unique lattice in $X'$ such that
$\cl\cap\cl'\sub\cl_1\sub\cl$, $\dim\cl/\cl_1=1$ and let $\cl_2$ be the unique lattice in $X'$ such that 
$\cl\cap\cl'\sub\cl_2\sub\cl'$, $\dim\cl'/\cl_2=1$. Let $\bar\cl_1=\cl_1/\e\cl_0$ (a line in $\cl/\e\cl$) and let
$\bar\cl_2=\cl_2/\e\cl'$ (a line in $\cl'/\e\cl'$). We say that $(\cl',g)\in Y'_n$ if $g(\bar\cl_1)=\bar\cl_2$; 
we say that $(\cl,g)\in Y''_n$ if $g(\bar\cl_1)\ne\bar\cl_2$. This defines the partition of $Y_n$. 

The subvarieties $Y_0,Y'_n (n>0),Y''_n (n>0)$ are precisely the pieces ${}^\tt\GG'/U_P$ of $Y=\GG'/U_P$. If $\kk$
is replaced by a finite field $F_q$ with $q$ elements then $Y_0,Y'_n,Y''_n$ become finite sets with cardinal 
$q(q^2-1),q^{2n}(q^2-1),q^{2n+1}(q^2-1)$ respectively.

\head 2. Character sheaves on $\GG^\d/U_P$\endhead
\subhead 2.1\endsubhead
We preserve the setup of 1.8. Let $K$ be a character sheaf on $C_\tt$ (the definition of \cite{\LIII} applies 
since $C_\tt$ is a connected component of an algebraic group whose identity component is $L_\tt$, a reductive 
algebraic group). Then $K$ is $(Q\cap Q')$-equivariant where $Q\cap Q'$ acts on $C_\tt$ through its quotient 
$L_\tt$ and there is a well defined $P$-equivariant simple perverse sheaf $K'$ on $P\T_{Q\cap Q'}C_\tt$ such that
$h_1^*K'=h_2^*K$ up to shift (here $P\T_{Q\cap Q'}C_\tt@<h_1<<P/(U_Q\cap U_{Q'})\T C_\tt@>h_2>>C_\tt$ are the
obvious maps). Now $K'$ can be viewed as a $P$-equivariant simple perverse sheaf on $\D$ (via $\c'$). Since 
$b:{}^\tt\GG^\d/U_P@>>>\D$ is a $P$-equivariant iterated affine space bundle we deduce that $b^*K'[d]$ is a 
$P$-equivariant simple perverse sheaf on ${}^\tt\GG^\d/U_P$ (for a well defined $d\in\ZZ$). Let 
$\ov{{}^\tt\GG^\d/U_P}$ be the closure of ${}^\tt\GG^\d/U_P$ in $\GG^\d/U_P$ (a $P$-stable algebraic subvariety of
$\GG^\d/U_P$). Let $\hK$ be the simple perverse sheaf on $\ov{{}^\tt\GG^\d/U_P}$ whose restriction to 
${}^\tt\GG^\d/U_P$ is $b^*K'[d]$. We can view $\hK$ as a simple perverse sheaf on $\GG^\d/U_P$ with support 
contained in $\ov{{}^\tt\GG^\d/U_P}$. 

Let $\cc_{P,\d}$ be the class of simple perverse sheaves on $\GG^\d/U_P$ consisting of all $\hK$ as above (for 
various $\tt$). The set of isomorphism classes in $\cc_{P,\d}$ is in bijection with the set of pairs $(\tt,K)$ 
where $\tt\in\ct(J,\d)$ and $K$ is a character sheaf on $C_\tt$ (up to isomorphism). The objects of $\cc_{P,\d}$ 
are called {\it character sheaves} on $\GG^\d/U_P$.

\subhead 2.2\endsubhead
Let $J\sub J'\sneq I$ and let $P\in\cp_J$, $Q\in\cp_{J'}$ be such that $P\sub Q$. Consider the diagram
$$\GG^\d/U_P@<f_1<<E_1@>f_2>>E_2@<f_3<<E_3@>f_4>>\GG^\d/U_Q$$
where 

$E_1=Q/U_Q\T\GG^\d/U_P$, $E_2=\{(P',gU_{P'})\in\cp_J\T\GG^\d/U_{P'};P'\sub Q\}$, 

$E_3=\{(P',gU_Q)\in\cp_J\T\GG^\d/U_Q;P'\sub Q\}$
\nl
are ind-varieties and 

$f_1(xU_Q,gU_P)=gU_P$, $f_2(xU_Q,gU_P)=(xPx\i,xgx\i U_{xPx\i})$, 

$f_3(P',gU_Q)=(P',gU_{P'})$, $f_4(P',gU_Q)=gU_Q$.
\nl
Note that $f_1$ is a principal fibration with group $Q/U_Q$; $f_2$ is a principal fibration with group $P/U_Q$;
$f_3$ is an affine space bundle with fibres isomorphic to $U_P/U_Q$; $f_4$ is a proper map.

Below we shall extend in an obvious way the standard operations for derived categories of $\bbq$-sheaves on
algebraic varieties to the case of ind-varieties.
Let $K\in\cc_{P,\d}$. Then $f_1^*K[d_1]=f_2^*K'[d_2]$ for a well defined simple perverse sheaf $K'$ on $E_2$. Here
$d_1=\dim Q/U_Q$, $d_2=\dim P/U_Q$. We set $\ind(K)=f_{4!}f_3^*K'$. By the decomposition theorem, $\ind(K)$ is a 
direct sum of simple perverse sheaves (with shifts) on $\GG^\d/U_Q$. More precisely, it is a direct sum of 
character sheaves (with shifts) on $\GG^\d/U_Q$. (This is shown by arguments similar to those in 
\cite{\LI, \S4,\S6}.) 

Conversely, for any character sheaf $\tK$ on $\GG^\d/U_Q$ there exists a character sheaf $K$ on $\GG^\d/U_P$ 
(where $P$ is an Iwahori subgroup contained in $Q$) such that some shift of $\tK$ is a direct summand of 
$\ind(K)$. (This is shown by arguments similar to those in \cite{\LI, \S4}.)


\widestnumber\key{\LI}
\Refs
\ref\key{\BE}\by R.B\'edard\paper On the Brauer liftings for modular representations\jour J.Algebra\vol93\yr1985
\pages332-353\endref
\ref\key{\LI}\by G.Lusztig\paper Parabolic character sheaves I\jour Mosc.Math.J.\vol4\yr2004\pages153-179\endref
\ref\key{\LII}\by G.Lusztig\paper Parabolic character sheaves II\jour Mosc.Math.J.\vol4\yr2004\pages869-896\endref
\ref\key{\LIII}\by G.Lusztig\paper Character sheaves on disconnected groups VI\jour Represent.Th.\vol8\yr2004
\pages377-413\endref
\ref\key{\LIV}\by G.Lusztig\paper Character sheaves on disconnected groups VIII\jour Represent.Th.\vol10\yr2006
\pages314-352\endref
\endRefs
\enddocument